\documentclass[12pt]{amsart}

\usepackage{latexsym}
\usepackage{amsmath}
\usepackage{amssymb}
\usepackage{amscd}
\usepackage{multicol}

\newtheorem{theo}{\bf Theorem}[subsection]
\newtheorem{theorem}{\bf Theorem}[section]
\newtheorem{corollary}[theorem]{\bf Corollary}
\newtheorem{coro}[theo]{\bf Corollary}%[subsection]

\newtheorem{propo}[theo]{\bf Proposition}%[subsection]
\newtheorem{lemma}[theo]{\bf Lemma}%[subsection]
\newtheorem{definition}[theo]{\bf Definition}%[subsection]
\newtheorem{definition-theorem}[theorem]{\bf Theorem-Definition}

\def\bR{\mathbb{R}}
\def\bZ{\mathbb{Z}}
\def\bC{\mathbb{C}}

\def\t{\mathfrak{t}}
\def\g{\mathfrak{g}}
 \def\k{\mathfrak{k}}

\def\g{\frak{g}}
\def\t{\frak{t}}
\def\p{\frak{p}}
\def\a{\frak{a}}

\def\u{\frak{u}}
\def\k{\frak{k}}

\def\o{\frak{o}}
\def\s{\frak{s}}

\numberwithin{equation}{section}

\title[Real loci of based loop groups]
{Real loci of based loop groups}

  \author[L. C. Jeffrey]{Lisa  C. Jeffrey}
   \address[L.\ C.\ Jeffrey]{Department of Mathematics\\ University of
 Toronto
    \\Toronto, Canada}
    \email{jeffrey@math.toronto.edu}

   \author[A.-L. Mare]{Augustin-Liviu  Mare}
  \address[A.-L. Mare]{Department of Mathematics and Statistics\\ University of Regina\\ Regina, Canada}
    \email{mareal@math.uregina.ca}

\begin{document}
\begin{abstract}
Let $(G,K)$ be a Riemannian symmetric pair of maximal rank, where
$G$ is a  compact simply connected Lie group and $K$ the fixed point set of an
involutive automorphism $\sigma$. 
This  induces an involutive automorphism $\tau$ of the based loop space
$\Omega(G)$.  There exists a maximal torus $T\subset G$ such that the
canonical action of $T\times S^1$ on $\Omega(G)$ is compatible with $\tau$
(in the sense of Duistermaat).
This allows us to formulate and  prove a  version of 
Duistermaat's   convexity theorem. Namely, the images of
$\Omega(G)$ and $\Omega(G)^\tau$ (fixed point set of $\tau$) under the 
$T\times S^1$ moment map on $\Omega(G)$ are
equal. The space $\Omega(G)^\tau$ is homotopy equivalent to the loop space $\Omega(G/K)$ of the Riemannian symmetric space
$G/K$.
We prove a stronger form of a result of Bott and Samelson which relates
the cohomology rings with coefficients in $\bZ_2$ of $\Omega(G)$ and $\Omega(G/K)$.
Namely, 
the two cohomology rings   are isomorphic, by a degree-halving isomorphism
(Bott and Samelson \cite{Bo-Sa} had proved that the Betti numbers are equal).
A version of this theorem involving equivariant cohomology is also proved.
The proof uses the notion of conjugation space in the sense of 
Hausmann,  Holm, and Puppe \cite{Ha-Ho-Pu}. 
\end{abstract}

\maketitle

\date{\today}

\tableofcontents

\newpage

\section{Introduction} 
Let $G$ be a compact connected simply connected Lie group.
Consider the space
$$\Omega(G):=\{\gamma : S^1 \to G \ : \ \gamma \ {\rm of \ Sobolev \ class \ } H^1,
\gamma(1)=e\}$$ of all based loops in $G$ (here $S^1$  is the unit circle in the complex plane). 
It is known that $\Omega(G)$ is an infinite dimensional symplectic manifold
which behaves in many respects like a {\it compact} symplectic manifold.
For example, let us consider the canonical action of the product $T\times S^1$
on $\Omega(G)$, %L02/18 
where
$T\subset G$ is a maximal torus and $S^1$ a circle (for more details, see 
Section \ref{duis}  %L02/18
below).  
One can show that this action is Hamiltonian. Moreover, by  the convexity theorem of Atiyah and Pressley \cite{At-Pr},
the image of the corresponding moment map 
 is a convex unbounded polyhedron
(by {\it convex polyhedron} we always mean  in this paper the convex hull of an infinite but discrete collection of points). %M 
Another instance of the same phenomenon is that the
$T\times S^1$-equivariant  cohomology of $\Omega(G)$ can be computed 
by Goresky-Kottwitz-MacPherson type %L02/18
formulas (this has been obtained by
Harada, Henriques, and Holm in \cite{Ha-He-Ho}).

Let $\sigma$ be a Lie group automorphism of $G$ with the following properties:
\begin{itemize}
\item $\sigma \circ \sigma ={\rm id}_G$, that is, $\sigma$ is an involution
\item there exists a maximal torus $T\subset G$ such that $\sigma(t)=t^{-1}$ for all
$t\in T$.
\end{itemize}
It is known (cf. e.g.  \cite[Chapter VI, Theorem 4.2]{Lo}) that 
any simply connected compact Lie group $G$ admits such an an automorphism 
$\sigma$. 
This $\sigma$ is unique up to an inner automorphism of  $G$.
For example, if $G=SU(n)$, $\sigma$ is given by 
$$\sigma((a_{k\ell})_{1\le k,\ell\le n})=(\bar{a}_{k\ell})_{1\le k,\ell\le n},$$
for any special unitary $n\times n$ matrix $(a_{k\ell})_{1\le k,\ell\le n}$
(the bar indicates the complex conjugate). Examples of such involutions 
for other Lie groups are presented in Section \ref{lasts}.

The automorphism $\sigma$ gives rise to the   involution $\tau$ of $\Omega(G)$ given by
\begin{equation}\label{taug}\tau(\gamma)(z)=\sigma(\gamma(\bar{z})),\end{equation}
for all $\gamma\in \Omega(G)$ and  $z \in S^1$. 
One can see that $\tau$ is an anti-symplectic automorphism of $\Omega(G)$,
that is, it satisfies $\tau^*(\omega)=-\omega$, where $\omega$ is the symplectic form 
of $\Omega(G)$ (cf. \cite{Ko}).
The automorphism $\tau$ of $\Omega(G)$ is compatible with the
$T\times S^1$ action mentioned above: that is, we have
\begin{equation}\label{co}\tau((t,z).\gamma) =(t^{-1},z^{-1}).\tau(\gamma),
\end{equation} for all $(t,z)\in T\times S^1$
and all $\gamma\in \Omega(G)$ (see Proposition \ref{propo} below). 
Real loci of compact (finite dimensional) symplectic manifolds with compatible
torus actions have been
investigated by several authors, like
Duistermaat \cite{Du}, O'Shea and Sjamaar \cite{OS-Sj},  %L02/18
Biss, Guillemin, and Holm \cite{Bi-Gu-Ho}, and Hausmann, Holm, and
Puppe \cite{Ha-Ho-Pu}.
 The loop space $\Omega(G)$ is
 infinite dimensional, thus we cannot  directly apply the results
 in the papers above. The goal of our paper is to  show that the following two results 
can be extended to $\Omega(G)$: the Duistermaat convexity theorem
(cf. \cite{Du}, see also Theorem \ref{dui} below)
and a more recent result of Hausmann, Holm, and Puppe which relates 
the (equivariant) %L02/18
 cohomology rings of the manifold and of the fixed point set of the involutive
automorphism
(cf. \cite{Ha-Ho-Pu}). More precisely, we  prove Theorems \ref{firstmain} and
\ref{secmain} below.
The first theorem concerns the moment map of the $T\times S^1$ action on
$\Omega(G)$, which is a map $\Omega(G)\to ({\rm Lie}(T)\oplus \bR)^*$.
 The explicit description of this map is given in Section \ref{convex}.
It turns out that it is more convenient to describe it
 by endowing ${\rm Lie}(G)$ with an ${\rm Ad}(G)$-invariant inner product
and restricting it to ${\rm Lie}(T)$, and then endowing
$\bR$ with the canonical inner product: we identify in this way
$ ({\rm Lie}(T)\oplus \bR)^*= {\rm Lie}(T)\oplus \bR$.

\begin{theorem}\label{firstmain} If $\Phi : \Omega(G)\to {\rm Lie}(T) \oplus \bR$ is the moment map
of the $T\times S^1$ action, then we have 
$$\Phi(\Omega(G))=\Phi(\Omega(G)^\tau).$$
Here $\Omega(G)^\tau$ denotes the fixed point set of $\tau$.
\end{theorem}

\noindent {\bf Remarks.}
1.  Let us consider the more general situation when $\sigma$ is an {\it arbitrary} involutive Lie group automorphism of $G$. The differential map $d\sigma_e$ is an involutive Lie algebra automorphism of 
$\g:={\rm Lie}(G)$. Let $\k,\p\subset \g$ denote the corresponding $+1$, respectively $-1$  eigenspaces. We have $\g=\k\oplus\p$. Let $\a$ be a maximal abelian subspace of $\g$ with $\a \subset \p$. Then $A:=\exp(\a)$ is a torus in $G$
(cf. e.g. \cite[Chapter VII]{He}). Let $T\subset G$ be a maximal torus such that
$A\subset T$. 
Consider again the involution $\tau$ of $\Omega(G)$ given by Equation (\ref{taug}).   
Again, $\tau$ is an antisymplectic automorphism of $\Omega(G)$.  The action of  $A\times S^1$ (which is a subgroup of $T\times S^1$) on $\Omega(G)$ is compatible with $\tau$. Let $\Phi_A : \Omega(G)\to \a\oplus \bR$ be the moment map of the
$A\times S^1$ action (as before, we make the identification 
$(\a\oplus \bR)^*=\a\oplus \bR$).  %M
It is not known whether
$\Phi_A(\Omega(G)^\tau)=\Phi_A(\Omega(G))$: this would be a version of
Duistermaat's convexity theorem stronger than Theorem \ref{firstmain} above.
Note that both $\Phi_A(\Omega(G))$ and $\Phi_A(\Omega(G)^\tau)$ are
convex polyhedra in $\a\oplus \bR$: the first by Atiyah and Pressley's theorem
mentioned above, the second by the convexity theorem of Terng \cite[Theorem 1.6]{Te-Convex} for infinite dimensional  
isoparametric submanifolds (for more details, see Section \ref{secterng} below).

 %M
2. It is probably also worth investigating whether the result in Theorem \ref{firstmain}
remains true if instead of  loops of Sobolev class $H^1$ we consider other classes of loops.
For example, let us consider the space $\Omega_{\rm alg}(G)$ of { algebraic loops}
in $G$ (see  Section \ref{equivc} for the exact definition of this notion). 
Note that $\Omega_{\rm alg}(G)$  is a $T\times S^1$ invariant subspace of $\Omega(G)$.  Atiyah and Pressley \cite{At-Pr} showed that 
we have $\Phi(\Omega(G))= \Phi(\Omega_{\rm alg}(G))$ (thus the latter set is 
also an unbounded convex polyhedron). The automorphism $\tau$ leaves
$\Omega_{\rm alg}(G)$ invariant.  We do not know whether 
$\Phi(\Omega_{\rm alg}(G)^\tau)=\Phi(\Omega_{\rm alg}(G))$. 
An important step towards the proof of this conjecture would be made by
taking the Bruhat cells in $\Omega_{\rm alg}(G)$ (see Section \ref{equivc} below) 
and their closures, which are finite dimensional projective varieties.
They are both $T\times S^1$ and $\tau$ invariant. One should first verify whether
for any such  variety $X$ we have $\Phi(X)=\Phi(X^\tau)$.

3.  The proof of Theorem 1.1 will be given in Section 2.
The main ingredients of the proof are as follows:
first, the set
$\Phi(\Omega(G)^\tau)$ is a convex subset of $\t\oplus \bR$
(as already mentioned in Remark 1 above, a slightly
more general result will be proved in Section 2.2);
second, the vertices of  Atiyah and Pressley's polyhedron
$\Phi(\Omega(G))$ are of the form $\Phi(\lambda)$, where
$\lambda : S^1 \to T$ is a group homomorphism
(see \cite[Section 1, Remark 2]{At-Pr}). 
 
 The following is the second main result of the paper.
 
\begin{theorem}\label{secmain} One has the following two ring isomorphisms:
\begin{align*}
{}&(a) \  H^{2*}(\Omega(G);\bZ_2)\simeq H^*
(\Omega(G)^\tau;\bZ_2)\\
{}& (b) \ H^{2*}_{T\times S^1}(\Omega(G);\bZ_2)\simeq H^*_{T_2\times \bZ_2}
(\Omega(G)^\tau;\bZ_2),
\end{align*}
where $T_2\times \bZ_2:=\{(t,z)\in T \times S^1 \ : \ t^2 =1 \ {\it  and} \ z=\pm 1\}$.
\end{theorem}

Note that the right-hand side of equation (b) above is well defined:  by the compatibility condition (\ref{co}), the group $T_2\times \bZ_2$ leaves
$\Omega(G)^\tau$ invariant. %L02/18

This theorem is 
related to a result of Bott and Samelson \cite{Bo-Sa} concerning the space of loops in a symmetric space.
To be more precise, let
 $G$ be, as before, a compact simply connected Lie group and 
 $\sigma$ a group automorphism of $G$ with the property that
 $\sigma\circ\sigma ={\rm id}_G$ (the assumption that $\sigma(t)=t^{-1}$ for all
 $t\in T$ is temporarily dropped).
 Then 
  $K=G^\sigma$ (the fixed point set of $\sigma$) is a connected closed
  subgroup of $G$ and
the homogeneous space $G/K$ has a natural structure of a Riemannian symmetric space. Explicit formulas for the $\bZ_2$ Betti numbers of the loop space 
$\Omega(G/K)$ are given in \cite[Corollary 3.10]{Bo-Sa}.
This result also gives the $\bZ_2$ Betti numbers of $\Omega(G)^\tau$,
since  the latter space  is homotopy equivalent to   
$\Omega(G/K)$ (see for instance Proposition \ref{equiv} below). 

Let us now reinforce the assumption that
$\sigma(t)=t^{-1}$ for all $t\in T$. Then $G/K$ is  a symmetric space of maximal rank
(that is, ${\rm rank} \ G/K={\rm rank} \ G$). 
Under this assumption, Bott and Samelson proved that
\begin{equation}\label{botsa}{\rm dim} \ H^{2q}(\Omega(G);\bZ_2)={\rm dim} \ H^q(\Omega(G/K);\bZ_2),\end{equation}
for all $q\ge 0$ (see \cite[Proposition 4.1]{Bo-Sa}).
The homotopy equivalence between    $\Omega(G)^\tau$ and $\Omega(G/K)$ mentioned above is  
$(T_2\times \bZ_2)$-equivariant with respect to a certain natural action of
$T_2\times \bZ_2$ on $\Omega(G/K)$ %L02/18
 (see Section \ref{two} below, especially 
Proposition \ref{equiv}). Consequently,
$\Omega(G)^\tau$ and $\Omega(G/K)$ have the same cohomology rings, both
equivariant and non-equivariant.
  In this way, the following result  can be deduced  
    from  Theorem  \ref{secmain}. Before stating it, we note that it is a stronger form of the result given by Equation (\ref{botsa}).
\begin{corollary}\label{bottsam}  If $G/K$ is a symmetric space of maximal rank, then
one has the following two ring isomorphisms:
\begin{align*}
{}&(a) \  H^{2*}(\Omega(G);\bZ_2)\simeq H^*
(\Omega(G/K);\bZ_2)\\
{}& (b) \ H^{2*}_{T\times S^1}(\Omega(G);\bZ_2)\simeq H^*_{T_2\times \bZ_2}
(\Omega(G/K);\bZ_2).
\end{align*}
\end{corollary}

\noindent {\bf Remarks.} 
1. The following result was also proved by Bott and Samelson. 
As usual, $G/K$ is a Riemannian symmetric space of maximal rank. Take $x \in \t$ and the orbits
${\rm Ad}_G(G)x=G/G_x$ and   ${\rm Ad}_G(K)x=K/K_x$.   
We have
$${\rm dim} \ H^{2q}(G/G_x;\bZ_2) = {\rm dim} \ H^q(K/K_x;\bZ_2),$$
for all $q\ge 0$ (see \cite[Proposition 4.3]{Bo-Sa}).
Stronger forms of this result have been obtained by Hausmann, Holm, and  Puppe in \cite{Ha-Ho-Pu}. Namely, they proved the
following ring isomorphisms:
\begin{align}
{}&\  H^{2*}(G/G_x;\bZ_2)\simeq H^*\label{unu}
(K/K_x;\bZ_2)\\
{}&  \ H^{2*}_{T}(G/G_x;\bZ_2)\simeq H^*_{T_2}
(K/K_x;\bZ_2)\label{doi},
\end{align}
where $T_2:=\{t\in T \ : \ t^2=1\}$.
The main idea of their proof is that $\sigma$ induces an anti-symplectic
involutive automorphism of  
$G/G_x$, which is compatible with the $T$ action and whose fixed point set
is $K/K_x$; the upshot is that this automorphism together with the
Schubert cell decomposition makes $G/G_x$ into a {\it spherical conjugation
complex}, and this automatically implies the isomorphisms (\ref{unu}) and (\ref{doi}).
Our proof of Theorem \ref{secmain} uses a similar argument.
Namely, we use the Bruhat cell decomposition of the space  $\Omega_{\rm alg}(G)$
in order to show that this, 
together with the involution $\tau$, is a spherical conjugation complex.
(We take this opportunity to note these arguments show that Theorem \ref{secmain} is also valid if $\Omega(G)$ is replaced by $\Omega_{\rm alg}(G)$.) 
Finally, we use a theorem which says that the inclusion $ \Omega_{\rm alg}(G)\hookrightarrow \Omega(G)$ is a 
homotopy equivalence. 
The details can be found in Section \ref{lastsec}.

%M
2. The following result can also be proved by using the methods of our paper,
combined with a theorem of Franz and Puppe (see  \cite[Theorem 1.3]{Fr-Pu}). 
Let $\kappa$ denote any of the two isomorphisms given at points (a) and (b) of
Corollary \ref{bottsam}, which are maps 
from the (equivariant) cohomology of $\Omega(G)$ to the (equivariant) cohomology
of $\Omega(G/K)$. Then we have
$$\kappa \circ {\rm Sq}^{2q} = {\rm Sq}^q \circ \kappa$$
for all $q\ge 0$. Here  $ {\rm Sq}^{2q}$ and ${\rm Sq}^q$
denote the Steenrod squaring operations on the (equivariant) cohomology rings
of $\Omega(G)$, respectively $\Omega(G/K)$.

\noindent {\bf Note.} By $S^1$ we will interchangeably denote the unit circle
in the complex plane and the quotient space $\bR/2\pi \bZ$. It will be clear from the context which of these two presentations is used.

\noindent {\bf Acknowledgements.} We would like to thank the referees for 
reading carefully previous versions of the manuscript and making several valuable suggestions.
 
 \section{The image of $\Omega(G)^\tau$ under the moment map}\label{convex}

\subsection{Duistermaat type convexity for 
$(\Omega(G),\tau,T\times S^1)$}\label{duis} %L02/18
Duistermaat proved the following theorem: 

\begin{theo}\label{dui}{\rm (\cite{Du})} Let $M$ be a compact symplectic manifold
 equipped with a Hamiltonian action of a torus ${\mathcal T}$ and an antisymplectic involution
 $\rho$   which are compatible, in the sense that
 \begin{equation}\label{compa} \rho(tx)=t^{-1}\rho(x),
 \end{equation} for all $t\in{\mathcal T}$ and all $x\in M$.
 If $\mu : M \to {\rm Lie}({\mathcal T})^*$ is the moment map of the ${\mathcal T}$  
  %L02/18
 action,
 then we have $$\mu(M)=\mu(M^\rho),$$
 where $M^\rho$ is the fixed point set of $\rho$.
\end{theo}

Our Theorem \ref{firstmain} is an extension of this result.
In this section we prove Theorem \ref{firstmain}. 
The considerations made in the introduction right before stating this theorem are
in force here.
We denote by $\g$ the Lie algebra of $G$ and   choose an ${\rm Ad}(G)$ invariant inner product on $\g$ (e.g. the negative of the
Killing form): if $X\in \g$ then $|X|$ denotes the length of $X$.

We consider  the action of $T$ on $\Omega(G)$ given by pointwise conjugation of loops,
that is,
\begin{equation}\label{conju}(t.\gamma)(\theta) =t\gamma(\theta) t^{-1},\end{equation}
for all $\gamma\in \Omega)G)$, $t\in T$,  %L02/15
and $\theta\in S^1$.
There is also an action of $S^1$ on $\Omega(G)$, given by the rotation of loops.
Concretely, if $e^{i\varphi}\in S^1$ and $\gamma\in \Omega(G)$, then
\begin{equation}\label{ei}(e^{i\varphi}.\gamma)(\theta):=\gamma(\theta+\varphi)\gamma(\varphi)^{-1}\end{equation}
for all $\theta\in S^1$. 

The details concerning the following results can be found for instance in
\cite{At-Pr}. First, the moment map of the $T$ action on $\Omega(G)$ is 
$p: \Omega(G) \to \t$ given by\footnote{The factors $\frac{1}{4\pi}$ in Equation (\ref{energy}) and
$\frac{1}{2\pi}$ in Equation (\ref{projection}) are due to a canonical choice of the
symplectic form on $\Omega(G)$, cf. e.g. \cite{At-Pr}.}
$$p(\gamma) = \frac{1}{2\pi}\int_0^{2\pi} {\rm pr}_\t
(\gamma(\theta)^{-1}\gamma'(\theta))d\theta
={\rm pr}_\t\left(\frac{1}{2\pi}\int_0^{2\pi} \gamma(\theta)^{-1}\gamma'(\theta)d\theta
\right)
,$$
where $\t$ is the Lie algebra of $T$ and ${\rm pr}_\t : \g \to \t$ is the orthogonal projection.
Second, the moment map of the $S^1$ action on $\Omega(G)$ is
the energy functional
$E:\Omega(G)\to \bR$,
\begin{equation}\label{energy}E(\gamma)=\frac{1}{4\pi}\int_0^{2\pi}|\gamma(\theta)^{-1}\gamma'(\theta)|^2
d\theta.\end{equation}
The actions of $T$ and $S^1$ commute with each other.
The moment map of the $T\times S^1$ action is
$$\Phi=p\times E : \Omega(G)\to \t \oplus \bR.$$
The following theorem was proved by Atiyah and Pressley  in \cite{At-Pr}:

\begin{theo}\label{atpr} {\rm (\cite{At-Pr})} We have
$$\Phi(\Omega(G))
={\rm cvx}\{\Phi(\lambda)  \ : \ 
\lambda :S^1 \to T \ {\it is \ a \ group \ homomorphism }  \}$$
where {\rm cvx} stands for convex hull. 
\end{theo}

We note that the group homomorphisms $S^1 \to T$ are precisely 
the elements of $\Omega(G)$ which are fixed by the $T\times S^1$  action. 

An important ingredient of this section  is the following result, which is a
consequence of the convexity theorem of Terng (see
\cite{Te-Convex}). We postpone its proof to  Section \ref{secterng} below.

\begin{theo}\label{terng}  The space $\Phi(\Omega(G)^\tau)$ is a convex subset of
$\t\oplus \bR$.
\end{theo}

%Now let us note that the space $\Phi_A(\Omega(G))$ is contained in
%the closure of the interior of the paraboloid\footnote{That is, in the region described by
 %$x \ge \frac{1}{2}||v||^2$.}
%$$\{(v,x)\in \a \times \bR \ : \ x = \frac{1}{2}||v||^2\}.$$
%Thus $\Phi_A(\Omega(G))$ is the convex hull of the 
%group homomorphisms $\gamma: S^1\to T$ such that
%$\Phi_A(\gamma)$ is on the boundary of the paraboloid, that is
%$$E(\gamma) =\frac{1}{2}||p_A(\gamma)||^2.$$

We are now ready to prove Theorem \ref{firstmain}.

\noindent{\it Proof of Theorem \ref{firstmain}.} First note that
$$\Phi(\Omega(G)^\tau)\subset \Phi(\Omega(G)).$$
To prove the opposite inclusion, 
we note that if $\lambda : S^1\to T$ is a group homomorphism,
then $\lambda \in \Omega(G)^\tau$. 
Indeed, for any $z \in S^1$ we have  
$$\tau(\lambda)(z)=\sigma(\lambda(\bar{z}))
=\sigma(\lambda(z^{-1}))=\sigma(\lambda(z)^{-1})
=\lambda(z).$$
From Theorem \ref{atpr} we deduce that $\Phi(\Omega(G))$ is the convex hull of some points
which are in 
$\Phi(\Omega(G)^\tau)$. 
Since the latter set is convex (by Theorem \ref{terng}), we deduce that
$\Phi(\Omega(G))\subset \Phi(\Omega(G)^\tau)$.
This finishes the proof.
\hfill $\square$

\subsection{Convexity  for $(\Omega(G)^\tau,A\times S^1)$}\label{secterng}%L02/18

The goal of this subsection is to prove Theorem \ref{terng}. 
In fact we will prove a stronger form of it. Namely, we consider
the situation described in  Remark 1 following Theorem \ref{firstmain}
and we show as follows:

\begin{theo}\label{terng2}  The space $\Phi_A(\Omega(G)^\tau)$ is a convex
subset of $\a\oplus \bR$.
\end{theo}

The following  expression of the moment map
$\Phi_A: \Omega(G)\to \a\oplus \bR$ will be needed in the proof (it can be deduced immediately from
the description of $\Phi: \Omega(G)\to \t\oplus \bR$ given in the previous subsection):
we have $\Phi_A=p_A\times E$, where
\begin{equation}\label{projection}p_A(\gamma) = \frac{1}{2\pi}\int_0^{2\pi} {\rm pr}_\a
(\gamma(\theta)^{-1}\gamma'(\theta))d\theta
={\rm pr}_\a\left(\frac{1}{2\pi}\int_0^{2\pi}\gamma(\theta)^{-1}\gamma'(\theta)d\theta
\right).\end{equation}

We also need the following considerations, which can be found
 in \cite{Te-PolarActions}. We  consider the loop group
$$L(G)=\{\gamma : S^1 \to G \  : \ \gamma \ {\rm of \ Sobolev \ class \ } H^1\}.$$
It acts by ``gauge transformations" on the Hilbert space $H^0(S^1,\g)$,  by
\begin{equation}\label{gam}\gamma \star u =\gamma u\gamma^{-1}-\gamma'\gamma^{-1},\end{equation}
for all $\gamma\in L(G)$ and $u\in  H^0(S^1,\g)$. The stabilizer  of the constant loop 
$0\in H^0(S^1,\g)$ consists of all $\gamma \in L(G)$ with $\gamma'\gamma^{-1}=0$,
which means that $\gamma$ is a constant loop in $G$. We deduce that the
$L(G)$ orbit
of $0$ can be identified with $L(G)/G$, which is the same as $\Omega(G)$.
Henceforth we will make the identification
\begin{equation}\label{ident}\Omega(G)=L(G)\star 0=\{\gamma'\gamma^{-1} \ : \ \gamma\in \Omega(G)\},\end{equation}
which is a subspace of $H^0(S^1,\g)$: more precisely, any based loop 
$\gamma:S^1\to G$ of Sobolev class $H^1$  is identified with 
$\gamma\star 0=\gamma^{-1}\gamma'$, which is an element  of $H^0(S^1,\g)$.
In this way, the moment map corresponding to the $T\times S^1$ action 
on $\Omega(G)$ is 
$\Phi : \Omega(G)\to \t\oplus \bR$, 
\begin{equation}\label{phiu}\Phi(u) =(P_\t(u), \frac{1}{2}\|u\|^2),\end{equation} for all
$u\in \Omega(G)$. Here  we regard $\t$ as a subspace of
$H^0(S^1,\g)$ (consisting of constant loops) and we denote
by $P_\t : H^0(S^1,\g)\to \t$  the orthogonal projection
with respect to the   canonical inner product on $H^0(S^1,\g)$.
We recall that this is given by
\begin{equation}\label{product}( u,v) =
\frac{1}{2\pi}\int_0^{2\pi}\langle u(\theta),v(\theta)\rangle d\theta,\end{equation}
for all $u,v\in H^0(S^1,\g)$ (here $\langle \ , \ \rangle$ is the ${\rm Ad}(G)$ invariant
inner product on $\g$  we chose at the beginning of this section). 
By $\| \cdot \|$ we denote the corresponding norm on $H^0(S^1,\g)$. %L02/18
To justify Equation (\ref{phiu}), we  show that
\begin{align}{}&P_\t(u)=\frac{1}{2\pi}\int_0^{2\pi}{\rm pr}_\t(u(\theta))d\theta\label{uno}\\
{}& \|u\|^2=\frac{1}{2\pi}\int_0^{2\pi}|u(\theta)|^2d\theta\label{due},
\end{align}
for all $u\in H^0(S^1,\g)$ (see also Equations (\ref{energy}) and (\ref{projection})).
Equation (\ref{due}) follows immediately from (\ref{product}).
To prove (\ref{uno}), we consider an orthonormal basis $e_1,\ldots, e_r$ of
$\t$, in the sense that $\langle e_i,e_j\rangle =\delta_{ij}$, for all
$1\le i,j\le r$ (here $\delta_{ij}$ is the Kronecker delta).  By using Equation
(\ref{product}), we deduce that
$(e_i,e_j)=\delta_{ij}$, for all $1\le i,j\le r$. Thus
\begin{align*}{}&P_\t(u)=\sum_{i=1}^r(u,e_i)e_i
= \sum_{i=1}^r\left(\frac{1}{2\pi}\int_0^{2\pi}\langle u(\theta),e_i\rangle d\theta\right) e_i\\{}&
 \ \ \ \ \ \ \ =\frac{1}{2\pi}\int_0^{2\pi}\sum_{i=1}^r\langle u(\theta),e_i\rangle e_i d\theta
=\frac{1}{2\pi}\int_0^{2\pi}{\rm pr}_\t(u(\theta))d\theta.
\end{align*}
Equation (\ref{phiu}) is now completely justified.

We recall now that $\sigma$ is an involution of $G$ whose fixed point set is $K$. We denote
$$\hat{K}:=\{\gamma\in L(G) \ : \ \gamma(-\theta)=\sigma(\gamma(\theta)), \ {\rm for\ all \ }
\theta\in S^1\}.$$
This is a subgroup of $L(G)$ which leaves invariant the closed vector subspace 
$$\hat{\p}(\g,\sigma):=\{u\in H^0(S^1,\g) \ : \ u(-\theta) =-d\sigma_e(u(\theta)), 
\ {\rm for\ all \ } \theta\in S^1\}$$
of $H^0(S^1,\g)$. 
As before, $\a$ is a maximal abelian subspace  of $\p$.  
It  can be made into a subspace of  $\hat{\p}(\g,\sigma)$ by regarding every element of
$\a$ as a constant loop.  In what follows we will need the notion of
{\it isoparametric submanifold} in Hilbert space.  %M
By definition, this is a finite codimensional Riemannian submanifold 
for which the normal vector bundle is flat relative to the normal connection and
satisfies some other assumptions: for instance, if $v$ is a parallel normal vector field
on the manifold, then the shape operators $A_{v(p)}$ and $A_{v(q)}$
corresponding to any two points $p$ and $q$ on the manifold are
 orthogonally conjugate.
For the exact definition  we refer the reader to 
 \cite[Section 6]{Te-Proper} (see also Chapter 7 of the monograph \cite{Pa-Te}).
We note that any isoparametric  submanifold induces a foliation of  Hilbert space by
parallel submanifolds\footnote{These are not necessarily isoparametric submanifolds.}, which we will call below the isoparametric foliation. 
%L02/18

\begin{propo}\label{isopara} 
(a) The orbits of the $\hat{K}$ action on $\hat{\p}(\g,\sigma)$ are elements of an isoparametric
foliation of the Hilbert space $\hat{\p}(\g,\sigma)$. 

(b) There exists $a\in \a$ such that the orbit 
$\hat{K}\star a$ is an isoparametric submanifold of $\hat{\p}(\g,\sigma)$.
The normal space at $a$ to this submanifold is $\a$.
\end{propo} 

\begin{proof} We use the following identifications (see also
\cite[Remark 3.4]{Te-Variational}):
\begin{align*}{}&\hat{K}=\{\gamma:[0,\pi]\to G \ : \ \gamma(0),\gamma(\pi)\in K\}
=:P(G,K\times K)\\
{}&\hat{\p}(\g,\sigma)=H^0([0,\pi], \g).
\end{align*}
By \cite[Theorem 1.2]{Te-PolarActions}, the action of $P(G,K\times K)$ on 
$H^0([0,\pi],\g)$ given by (\ref{gam}) is polar
(by definition, which can be found in full detail in \cite{Te-PolarActions}, this means essentially that there exists a section of this action, that is, a submanifold of
$H^0([0,\pi],\g)$ which meets all  orbits of the action and meets them orthogonally). 
%L02/18
By \cite[Theorem 8.10]{Te-Proper}, the orbits
of this action are an isoparametric foliation. In particular, the principal orbits are isoparametric submanifolds. We are looking for such orbits.
To find them, we recall  that the action of $K\times K$ on $G$ given by
$$(k_1,k_2).g=k_1gk_2^{-1},$$
for all $k_1,k_2\in K$ and $g\in G$ is polar; a section of this action is $A=\exp(\a)$ (cf. e.g. \cite{Co}). 
From \cite[Theorem 1.2]{Te-PolarActions} we deduce that $\a$ (the space
of constant maps from $[0,\pi]$ to $\a$) is a section
of the $P(G,K\times K)$ action on $H^0([0,\pi],\g)$.
To prove our proposition, we only need to pick $a\in \a$ a regular point (that is, one whose orbit is principal). Such a point exists due to the following
criterion (see \cite[Theorem 1.2, (6)]{Te-PolarActions}): a point $a\in\a$ is regular for the $P(G,K\times K)$ action on $H^0([0,\pi],\g)$
if and only if $\exp(a)$ is regular for the $K\times K$ action on $G$. Moreover, a general result says that any section 
of a polar action of a compact Lie group on a simply connected compact manifold contains
regular points (see e.g. \cite[Theorem 1.6]{Te-PolarActions}).
This finishes the proof.
\end{proof}

In order to  prove Theorem \ref{terng} we will show that, via the identification
(\ref{ident}), $\Omega(G)^\tau$ is the same as the element $\hat{K}\star 0$ of the
isoparametric foliation in the previous proposition. Then we use the convexity theorem
for isoparametric foliations of Terng \cite{Te-Convex}.
For the moment, we will prove the following lemma. 

\begin{lemma} Take $\gamma\in \Omega(G)$ and denote
$\gamma_0=\tau(\gamma)$. Then we have
$${\gamma_0}'(\theta){\gamma_0}^{-1}(\theta)=
-d\sigma_e(\gamma'(-\theta)\gamma^{-1}(-\theta))$$
for all $\theta\in S^1$.
\end{lemma}

\begin{proof} If $g\in G$, then the tangent space to $G$ at $g$
consists of vectors of the form $Xg=dR(g)_e(X)$, where $X\in T_eG$.
Here $R(g):G\to G$ is the right multiplication by $g$.
Moreover, we have
$$d\sigma_g(Xg)=d\sigma_e(X)\sigma(g).$$
Indeed,
$$d\sigma_g(Xg)=d\sigma_g(dR(g)_e(X))
=d(\sigma \circ R(g))_e(X)
=d(R(\sigma(g))\circ \sigma)_e(X)
=d\sigma_e(X)\sigma(g).$$
We deduce that
\begin{align*}{}&{\gamma_0}'(\theta){\gamma_0}^{-1}(\theta)
=d\sigma_{\gamma(-\theta)}(-\gamma'(-\theta))
\sigma(\gamma(-\theta)^{-1})\\
{}&
=-d\sigma_{\gamma(-\theta)}(\gamma'(-\theta)\gamma(-\theta)^{-1}\gamma(-\theta))
\sigma(\gamma(-\theta)^{-1})
=-d\sigma_{e}(\gamma'(-\theta)\gamma(-\theta)^{-1}).\end{align*}
\end{proof}

From this lemma we deduce
\begin{equation}\label{o}\Omega(G)^\tau=\{u\in \Omega(G) \ : \ -d\sigma_e(u(\theta))=u(-\theta) \ {\rm for  \ all \ } 
\theta\in S^1\}
=\Omega(G)\cap\hat{\p}(\g,\sigma).\end{equation}
This space is the same as the orbit $\hat{K}\star 0$, as the following lemma shows.

\begin{lemma}\label{fixq} We have
$$\Omega(G)^\tau=\hat{K}\star 0.$$
\end{lemma}

\begin{proof} The inclusion $\hat{K}\star 0\subset \Omega(G)^\tau$ is clear,
because $\hat{K}\star 0$ is a subset of both $\hat{\p}(\g,\sigma)$ and
$L(G)\star 0$.  We now prove the reverse inclusion. 
Take $\gamma\in \Omega(G)^\tau$: by identifying it 
with the element $\gamma\star 0=\gamma'\gamma^{-1}$ of $H^0(S^1,\g)$ and taking into account
Equation (\ref{o}), we have 
$$d\sigma_e(\gamma'(\theta)\gamma^{-1}(\theta))
=-\gamma'(-\theta)\gamma^{-1}(-\theta),$$
for all $\theta\in S^1$. We show that $\gamma\in \hat{K}$,
as follows.
We have
\begin{align*}{}&\frac{d}{d\theta}\sigma(\gamma(\theta))
\\=&d\sigma_{\gamma(\theta)}(\gamma'(\theta))
\\=&d\sigma_e(\gamma'(\theta)\gamma^{-1}(\theta))\sigma(\gamma(\theta))
\\=&-\gamma'(-\theta)\gamma^{-1}(-\theta)\sigma(\gamma(\theta)),
\end{align*}
which implies
$$\frac{d}{d\theta}(\sigma(\gamma(\theta)))\sigma(\gamma(\theta))^{-1}
=\frac{d}{d\theta}(\gamma(-\theta))\gamma(-\theta)^{-1}.$$
We deduce that the loops $\theta \mapsto \sigma(\gamma(\theta))$
and $\theta\mapsto \gamma(-\theta)$ are equal. Thus
$\tau(\gamma)=\gamma$, in other words, $\gamma \in \hat{K}$.
\end{proof}

We are now ready to prove our main result.

\noindent {\it Proof of Theorem \ref{terng2}.} By the convexity theorem of
Terng (see
\cite[Theorem 1.6]{Te-Convex}),
the image of the map
$\label{mapp}\Psi_A:\hat{K}\star 0 \to \a \oplus \bR$ given by
$$ \Psi_A(u)= (P_\a(u), \|u\|^2)
$$
is a convex polyhedron  in
$\a\oplus \bR$ (we are also using Proposition \ref{isopara}). 
Here $P_\a : \hat{\p}(\g,\sigma)\to \a$ is the orthogonal projection with respect to the
Hilbert space metric.
By Lemma \ref{fixq}, $\Psi_A(\Omega(G)^\tau)$ is a convex polyhedron.  
If we now compare the map $\Psi_A$ with the moment map
$\Phi_A=p_A \times E$ (see Equations (\ref{energy}) and (\ref{projection})),
we note that the two maps are essentially the same. More specifically,
by taking into account the identification given by (\ref{ident}), we have
$$\Phi_A(u)=\left(P_\a(u), \frac{1}{2}\|u\|^2\right),
$$
for all $u\in \Omega(G)^\tau$ (this can be proved in the same way as
equation (\ref{phiu})). 
We deduce  that the set $\Phi_A(\Omega(G)^\tau)$ is obtained
from $\Psi_A(\Omega(G)^\tau)$ by the  automorphism
of $\a\oplus \bR$ given by
$$(a,r)\mapsto \left(a, \frac{1}{2}r\right)$$
for all $(a,r)\in \a\oplus \bR$. Thus $\Phi_A(\Omega(G)^\tau)$ is a
convex polyhedron as well. This finishes the proof.  %L02/18
\hfill $\square$

\section{(Equivariant) cohomology ring of $\Omega(G/K)$}\label{lastsec}

\subsection{(Equivariant) cohomology of $\Omega(G)^\tau$}\label{equivc}
In this subsection we will prove Theorem \ref{secmain}. An important ingredient of the proof
will be the space $\Omega_{\rm alg}(G)$ of algebraic loops in $G$. 
By definition, this is
$$\Omega_{\rm alg}(G)=L_{\rm alg}(G^{\bC})\cap \Omega(G).$$
Here $G^{\bC}$ denotes the complexification of the Lie group $G$
and $L_{\rm alg}(G^{\bC})$ is the set of all (free)
loops $\gamma : S^1\to G^\bC$ which are restrictions  of algebraic maps  
from $\bC^*$ to $G^{\bC}$. In the case when $G^{\bC}$ is a subgroup
of some general linear group $GL_n(\bC)$, elements of $L_{\rm alg}(G^{\bC})$ are
Laurent series of the form 
\begin{equation}\label{degree}\gamma(z)=\sum_{p=-k}^k z^pA_p,\end{equation}
for some $k\ge 0$, where $A_p$ are elements of ${\rm Mat}^{n\times n}(\bC)$. 
For a fixed $k$, the space of all maps $\gamma$ of the form (\ref{degree}) is equipped with the standard metric topology which comes from its identification
with $\left({\rm Mat}^{n\times n}(\bC)\right)^{2k+1}$; we denote by
$\Omega_{\rm alg}^k(G)$ the space of all $\gamma$ of type (\ref{degree}) which map $S^1$ to $G$,
and equip it with the subspace topology.  
We endow $\Omega_{\rm alg}(G)$ with the direct limit topology coming from the
filtration $\{\Omega_{\rm alg}^k(G)\}_{k\ge 0}$.
The following theorem has been proved by Mitchell in
\cite{Mi} (see
Theorem 4.1 and the theorem in the introduction of his paper, where the result is attributed to
Quillen).   Another proof can be found in  \cite{Ko} (see Theorem 3.1.4 of that paper).

\begin{theo}\label{mitc}{\rm (\cite{Mi}, \cite{Ko})} (a) The inclusion map $\Omega_{\rm alg}(G) \to \Omega(G)$ is
a homotopy equivalence.

(b) The automorphism $\tau$ of $\Omega(G)$ leaves
 $\Omega_{\rm alg}(G)$ invariant
and the inclusion map 
$\Omega_{\rm alg}(G)^\tau \to \Omega(G)^\tau$ is
a homotopy equivalence.
\end{theo}

The advantage of dealing with $\Omega_{\rm alg}(G)$ instead of $\Omega(G)$ is that the
former space has a natural CW-decomposition. Its elements are the Bruhat  cells, which 
are described in what follows
(the details of this construction can be found in \cite[Sections 2 and 3]{Mi}). First we make the identification 
\begin{equation}\label{ome}\Omega_{\rm alg}(G)=L_{\rm alg}(G^{\bC})/L^+_{\rm alg}(G^{\bC})\end{equation}
 where $L^+_{\rm alg}(G^{\bC})$ 
is the subgroup of $L_{\rm alg}(G^{\bC})$ consisting of loops of the
form (\ref{degree}) for some $k\ge 0$, where $A_p=0$ for all $p< 0$.
We consider the roots of $G$ with respect to $T$, which are linear functions
$\t \to \bR$. 
The root space decomposition of $\g^{\bC}:=\g\otimes \bC$ is
$$\g^{\bC} = \t^\bC \oplus \sum_\alpha \g^{\bC}_\alpha,$$
where the sum runs over all the roots of $G$ with respect to $T$. 
We fix a simple root system $\alpha_1,\ldots,\alpha_\ell$ and
denote by $B^-$ the (Borel) connected subgroup of $G^{\bC}$ whose Lie algebra is 
$\t^ \bC \oplus \sum _\alpha \g^{\bC}_{\alpha}$, where 
the sum runs over 
all negative roots $\alpha$.
The Bruhat decomposition of $\Omega_{\rm alg}(G)$ is
\begin{equation}\label{cw}\Omega_{\rm alg}(G)=\bigsqcup_\lambda {\mathcal B}\lambda
\end{equation}
where the union runs over all group homomorphisms $\lambda: S^1\to T$
such that $\lambda'(0)$ is in the closure of the fundamental Weyl chamber of $\t$.
Here ${\mathcal B}$ is the subgroup of $L^+_{\rm alg}(G^{\bC})$ consisting of all loops 
$\gamma$ of the form (\ref{degree}) 
for some $k\ge 0$, where $A_p=0$ for all $p< 0$ and $A_0 \in B^-$.
The decomposition described by (\ref{cw}) is a CW decomposition
(cf. e.g. \cite[Section 3]{Mi}). The orbits ${\mathcal B}\lambda$ are the  Bruhat cells.

%M
Any Bruhat cell is
homeomorphic to some  complex vector space.
Proposition \ref{mitc} below gives a more precise description of this
homeomorphism.
In order to state it, we need to make some more considerations. 
First we note that the set of  group homomorphisms $\lambda: S^1 \to T$ 
can be identified with the integer  lattice  $I=\ker(\exp : \t \to T)$.
Let $W$ be the Weyl group of $G$. We recall that this is the group of linear transformations of
$\t$ generated by the
reflections about the hyperplanes $\ker\alpha_1$, $\ker\alpha_2,\ldots,\ker\alpha_\ell$;
let us denote these reflections by $s_1,s_2, \ldots,s_\ell$.
The affine Weyl group $\tilde{W}$ is the semidirect product 
$W\ltimes I$. It  is the same as the group of affine transformations of
$\t$  generated by 
$s_1,s_2,\ldots,s_\ell$, and $s_0$. Here $s_0$ is the reflection about
the affine hyperplane $\{x\in \t \ : \ \alpha_0(x)=1\}$, where
$\alpha_0$ is the highest root of $G$.  
To any $s\in \{s_0,s_1,\ldots,s_\ell\}$ we assign  the subgroup $U_s$
of $L_{\rm alg}(G^{\bC})$, as follows:
\begin{itemize}
\item For $j\in \{1,\ldots, \ell\}$ we have $U_{s_j}:=\exp (\g^{\bC}_{\alpha_j})$ 
(its elements are constant loops in $G^{\bC}$). Since $U_{s_j}$ is a unipotent 
group, the exponential map is an
isomorphism between $U_{s_j}$ and its Lie algebra  $\g^{\bC}_{\alpha_j}$. 
More precisely, by fixing $E_{\alpha_j}$ a non-zero vector in
$\g^{\bC}_{\alpha_j}$, the map $\bC\to U_{\alpha_j}$ given by 
$x\mapsto \exp(x E_{\alpha_j})$ is a homeomorphism.
\item $U_{s_0}$ consists of loops of the form $z\mapsto \exp(z^{-1} X)$, $z\in S^1$,
where $X\in \g^{\bC}_{-\alpha_0}$. Again, since $U_{s_0}$ is a unipotent group, the exponential map is an
isomorphism between $U_{s_0}$ and  $\g^{\bC}_{-\alpha_0}$.
By fixing again $E_{-\alpha_0}$ a non-zero vector in
$\g^{\bC}_{-\alpha_0}$, the map $\bC\to U_{\alpha_0}$ which assigns to
$x\in \bC$ the loop $z\mapsto \exp (z^{-1}xE_{-\alpha_0})$ is a homeomorphism.
\end{itemize}
We mention without any further explanations that the groups $U_s$ are the root
subgroups of $L_{\rm alg}(G^{\bC})$ corresponding to a certain canonical simple affine   root system of $G$ (note that the Lie algebra of $L_{\rm alg}(G^{\bC})$ has a root decomposition
labeled by the affine roots).   

%The affine Weyl group can also be realized as $$\tilde{W}=
%W\ltimes {\rm Hom}(S^1,T)=L_{\rm alg}(N^{\bC})/T^{\bC}$$ where
%$ {\rm Hom}(S^1,T)$ is the set of all group homomorphisms from $S^1$ to $T$,
%$N$  the normalizer of $T$ in $G$ and $L_{\rm alg}(N^{\bC})$ the set of algebraic maps from $\bC^*$ to the complexified Lie group $N^{\bC}$. 
%To any $k\in    \{0,1,\ldots,\ell\}$ we  assign an element
%$\tilde{n}_k \in L_{\rm alg}(N^{\bC})$ such that 
%$s_k=\tilde{n}_kT^{\bC}$, as follows:
%\begin{itemize}
%\item For $k\in \{1,\ldots,\ell\}$ we take $n_k\in N$ such that
%the element $s_k$ of $W=N/T$ is the same as the coset $n_kT$.
%We set $\tilde{n}_k:=n_k$.
%\item For $k=0$, we  take again $n_0\in N$ such that
%$s_0=n_0T$. This time we set
%$$\tilde{n}_0:=n_0\lambda_0,$$
%where $\lambda_0: S^1\to T$ is given by 
%$\lambda_0(\theta)=\exp(-\frac{\theta}{2\pi}\alpha_0^{\vee}),$ for all $\theta \in S^1=\bR/2\pi \bZ$.
%\end{itemize} 

Take $\lambda \in I =\tilde{W}/W$ and consider the element  $\tilde{w}$  of  $\tilde{W}$ 
which has minimal length (with respect to the generating set $s_0,s_1,\ldots,s_\ell$) and satisfies $\lambda = \tilde{w}W$. Let $\tilde{w}=s_{i_1}\ldots s_{i_k}$ be any reduced decomposition of $\tilde{w}$, where $i_1,\ldots,i_k\in \{0,1,\ldots,\ell\}$.
The following result has been proved by Mitchell in \cite{Mi}:
\begin{propo}{\rm (\cite{Mi})}\label{mitch}  The map 
\begin{align}{}\label{iden}{}&\bC^k=U_{s_{i_1}}\times \ldots \times U_{s_{i_k}}\to 
L_{\rm alg}(G^{\bC})/L^+_{\rm alg}(G^{\bC})=\Omega_{\rm alg}(G)\\{}&
(u_1,\ldots,u_k)\mapsto u_1 \ldots u_k L^+_{\rm alg}(G^{\bC})\nonumber
\end{align} is a homeomorphism
onto the Bruhat cell ${\mathcal B}\lambda$. 
\end{propo} 
Let $\sigma$ be the automorphism of $G$ defined in the introduction.
We note that the involutive automorphism $\tau$ of $\Omega(G)$ given by
(\ref{taug}) leaves $\Omega_{\rm alg}(G)$ invariant. 
To understand this, we first extend $\sigma$ to a group automorphism of
$G^{\bC}$, namely the one whose differential at the identity element is 
the anti-complex linear extension of the differential of the original $\sigma$. 
That is, we have
\begin{equation}\label{conj}d\sigma_e(X+iY)=d\sigma_e(X)-id\sigma_e(Y),\end{equation}
for all $X,Y\in \g$. 
Then we extend
$\tau$ to a group automorphism of $L_{\rm alg}(G^{\bC})$,
namely the one described by  
Equation (\ref{taug}) with  $\gamma$  in $L_{\rm alg}(G^{\bC})$.
This map  leaves $L^+_{\rm alg}(G^{\bC})$ invariant and
induces the original automorphism $\tau$ of $\Omega_{\rm alg}(G)$ via the identification
(\ref{ome}). 
 
We now consider the decomposition $\g=\k\oplus \p$, where
$\k=\{X\in \g \ : \ d\sigma_e(X)=X\}$ and $\p=\{X\in \g \ : \ d\sigma_e(X)=-X\}$.
Note that $\t$ is a subset of $\p$.
The automorphism $d\sigma_e$ of $\g^{\bC}$ has fixed point set equal to
$\g_0:=\k+i\p$. The latter space is a real form of $\g^\bC$. 
Any root $\alpha$ of $\g^\bC$ with respect to $\t\otimes \bC$
takes real values on the subspace $i\t$ of   $\g_0$. This means that $\g_0$ is a 
split real form of $\g^{\bC}$ (cf. e.g \cite[Section 26.1]{Fu-Ha}). We deduce that
we have the splitting
$$\g_0=i\t \oplus \sum \bR E_\alpha,$$
where the sum runs over all the roots $\alpha$ of $G$ with respect to
$T$ and  $E_\alpha$ is a (nonzero) root vector
for any such root $\alpha$.  In constructing the groups $U_{s_0},U_{s_1},\ldots,
U_{s_\ell}$ (see above) we  use the vectors $E_{-\alpha_0}, E_{\alpha_1},\ldots,
E_{\alpha_\ell}$ in the previous equation.

We will prove the following result
(see also \cite[Proof of Theorem 5.9]{Mi}).
 
 \begin{propo}\label{pr} Any Bruhat cell 
 ${\mathcal B}\lambda$  in $\Omega_{\rm alg}(G)$ remains invariant under  $\tau$.
 Moreover, via the homeomorphism  $\bC^k\simeq {\mathcal B}\lambda$ described by
 equation (\ref{iden}),
 $\tau$ acts on ${\mathcal B}\lambda$  by complex conjugation.
 \end{propo}
 
 \begin{proof} We have already seen that if $\lambda: S^1\to T$ is a group homomorphism then 
 $\tau(\lambda)=\lambda$ (see the proof of Theorem \ref{firstmain} at the end of Section
 \ref{duis}).
The automorphism $\tau$ leaves ${\mathcal B}$ invariant:
this follows from the definition
of ${\mathcal B}$ and the fact that the Borel subgroup $B^-$ is
$\sigma$-invariant. Consequently, $\tau$ leaves the orbit ${\mathcal B}\lambda$ 
invariant. The homeomorphism
$  U_{s_{i_1}}\times \ldots \times U_{s_{i_k}}\to
{\mathcal B}\lambda$ described by Equation (\ref{iden}) is $\tau$-equivariant, where 
$\tau$ acts diagonally on the
 domain of the map. The reason is that $\tau$ is a group automorphism of 
 $L_{\rm alg}(G^{\bC})$. 
The last statement in the proposition follows from the fact that
$\tau$ leaves  $U_{s_j}$ invariant, for any $j\in \{0,1,\ldots,\ell\}$;
moreover, via the identification $U_{s_j}=\g^{\bC}_{\alpha_j}=\bC$ (see above), $\tau$ acts as the complex conjugation. Indeed, if $j\neq 0$
then $\g^{\bC}_{\alpha_j}=\bC E_{\alpha_j}$ and by Equation 
(\ref{conj}), for any $x\in \bC$ we have
$$\tau(\exp(xE_{\alpha_j}))=\sigma(\exp(xE_{\alpha_j}))=
\exp(d\sigma_e(xE_{\alpha_j}))=
\exp(\overline{x}E_{\alpha_j});$$
 for $j=0$, we use that for any complex number
$x$, the loop $z\mapsto \exp (z^{-1}xE_{-\alpha_0})$ is mapped by $\tau$ to
$$z\mapsto \sigma(\exp(zxE_{-\alpha_0}))=\exp(\overline{zx}E_{-\alpha_0}))
=\exp(z^{-1}\overline{x}E_{-\alpha_0}).$$
\end{proof}

Our proof of Theorem \ref{secmain} uses the notion of 
spherical conjugation complex, defined in \cite{Ha-Ho-Pu}.
By definition, a spherical conjugation complex is a (finite or infinite) cell complex $X$  equipped with an involutive automorphism
$\rho$
with the following properties:
\begin{itemize}
\item each cell in $X$ is a complex cell, that is, it is homeomorphic to 
$\bC^k$, for some $k \in \bZ$, $k\ge 0$
\item $\rho$ leaves each cell $\bC^k$ invariant, acting on it as the
complex conjugation. That is, we have
$$\rho(z_1,\ldots, z_k) =(\bar{z}_1, \ldots, \bar{z}_k),$$
for all $(z_1,\ldots, z_k) \in \bC^k$.
\end{itemize}
The following theorem has been proved in \cite[Sections 5 and 7]{Ha-Ho-Pu}.

\begin{theo}\label{hhp} {\rm (\cite{Ha-Ho-Pu})} Let  $(X,\rho)$ be a spherical conjugation complex and denote by $X^\rho$  the fixed point set of $\rho$.
Then we have as follows:

(a) There exists a 
degree-halving ring isomorphism
$H^{2*}(X;\bZ_2) \simeq H^*(X^\rho;\bZ_2)$. 

(b) 
Let ${\mathcal T}$ be a compact torus acting on $X$ such that the action
is   compatible with $\rho$, in the sense that
$$\rho(tx)=t^{-1}\rho(x)$$
for all $t\in {\mathcal T}$ and all $x\in X$.
Then there exists a degree-halving ring isomorphism
$H^{2*}_{{\mathcal T}}(X;\bZ_2)\simeq H^{*}_{{\mathcal T}_2}(X^\rho;\bZ_2)$.
Here ${\mathcal T}_2$ denotes the set of all $t \in {\mathcal T}$ with $t^2=1$.  %L02/18
\end{theo} 
 
 Without any further comments we mention that the key point of this theorem is that a spherical conjugation complex is a
 conjugation space (for the definition of this notion, see \cite{Ha-Ho-Pu}).
 
 We are now ready to give the desired proof:
 
 \noindent {\it Proof of Theorem \ref{secmain}.} 
By
Proposition \ref{pr},   $\Omega_{\rm alg}(G)$ together with the
 involution $\tau$ is a spherical conjugation complex.  
 Theorem \ref{hhp} (a) implies that  we have a ring isomorphism
 $$H^{2*}(\Omega_{\rm alg}(G)) \simeq H^*(\Omega_{\rm alg}(G)^\tau).$$
 Combined with Theorem \ref{mitc}, this implies point (a) of
 Theorem \ref{secmain}.
 Point (b) follows from the fact that the actions of $T\times S^1$ and $\tau$ 
 on $\Omega_{\rm alg}(G)$ are compatible,
 see Proposition \ref{propo} below. We use Theorem \ref{hhp} (b) and again
 Theorem \ref{mitc}.
\hfill $\square$
 
\subsection{The Bott-Samelson theorem for $\Omega(G/K)$}\label{two}
Throughout this subsection 
$G$ will be a compact connected simply connected Lie group and
$\sigma$ an arbitrary involutive automorphism of $G$. The notations established
in Remark 1 following Theorem \ref{firstmain} are in force.
We consider again the group $K=G^\sigma$ and the homogeneous space $G/K$,
which has a canonical structure of a Riemannian symmetric space. 
Let us also consider the loop space\footnote{The reason why the loops in this definition are
defined on $[0,\pi]$, and not on $[0,2\pi]$ or $S^1=\bR/2\pi \bZ$, as usual,
will be understood later (see the proof of  %L02/18
Proposition \ref{equiv}).}
$$\Omega(G/K):=\{\mu : [0,\pi] \to G/K \ : \ \mu {\rm \ is \ of \ Sobolev \ class \ }H^1
\ {\rm and \ } \mu(0)=\mu(\pi)=eK\} $$
where $eK$ denotes the coset of $e$ in $G/K$.
Consider the group
$$A_2:=\{a\in A \ : \ a^2=e\}.$$
In this subsection we will define  $A_2\times \bZ_2$ actions on $\Omega(G/K)$
and $\Omega(G)^\tau$, show that these two spaces are equivariantly homotopy equivalent, and finally prove  Corollary \ref{bottsam}. 

  We first note that
$A_2=A\cap K.$
This can be justified as follows:
if $a\in A$ then  $a=\exp (X)$ where
$X\in \a$, so $\sigma(a)=a^{-1}$; consequently
$$a\in K \Leftrightarrow \sigma(a)=a 
 \Leftrightarrow a^{-1}=a  \Leftrightarrow a^2=e.$$
The group $A_2$ acts on $\Omega(G/K)$ by pointwise multiplication of
the loops from the left:
$$(a.\mu)(\theta) = a\mu(\theta),$$
for all $a\in A_2$, $\mu\in \Omega(G/K)$ and $\theta \in [0,\pi]$.
There is also an action of $\bZ_2$ on $\Omega(G/K)$, which is
more subtle. It is determined by the involutive automorphism
$\mu\mapsto \tilde{\mu}$ of $\Omega(G/K)$, defined below.
We first prove a lemma:
\begin{lemma}\label{lemmafirst}Any loop $\mu \in \Omega(G/K)$ can be written as
$$\mu(\theta) = \gamma(\theta)K,$$
where $\gamma : [0,\pi]\to G$ is an $H^1$ map such that
$\gamma(0)=e$ and $\gamma(\pi)\in K$.
\end{lemma}
\begin{proof} 
We use the Path Lifting Theorem (cf. e.g. \cite[Theorem 3.4.30]{Ab-Ma-Ra})
for the locally trivial bundle $G\to G/K$. 
\end{proof}
 
\begin{definition}\label{dfirst} Let $\mu\in \Omega(G/K)$ be of the form
$\mu(\theta)=\gamma(\theta)K$, $\theta\in [0,\pi]$, like
 in the previous lemma.  We define $\tilde{\mu}$ by
$$\tilde{\mu}(\theta):=\sigma(\gamma(\pi-\theta))K,$$
$\theta\in [0,\pi]$.
\end{definition}

We first verify that the map $\mu\mapsto \tilde{\mu}$ is independent of the
choice of $\gamma$: if $\gamma_1$ is another representative of $\mu$,
that is, if $\gamma_1(\theta)=\gamma(\theta)k$, for some $k\in K$,
then
$$\sigma(\gamma_1(\pi-\theta))K=\sigma(\gamma(\pi-\theta)k)K
=\sigma(\gamma(\pi-\theta))\sigma(k)K
=\sigma(\gamma(\pi-\theta))kK
=\sigma(\gamma(\pi-\theta))K.$$
Next we verify that  the map  $\mu\mapsto \tilde{\mu}$ is
involutive, that is $\tilde{\tilde{\mu}}=\mu$.
To do this, we write
$$\tilde{\mu}(\theta):=\sigma(\gamma(\pi-\theta))\gamma(\pi)^{-1}K,$$
and deduce that
$$\tilde{\tilde{\mu}}(\theta):=\sigma(\sigma(\gamma(\pi-(\pi-\theta))\gamma(\pi)^{-1}))K
=\gamma(\theta)K=\mu(\theta).$$
In this way we have defined our $\bZ_2$ action on $\Omega(G/K)$.

\begin{lemma}\label{commute} The $A_2$ and $\bZ_2$ actions on $\Omega(G/K)$ defined
above commute with each other and thus define an action of
$A_2\times \bZ_2$.
\end{lemma} 

\begin{proof} Take $a\in A_2$ and $\mu\in \Omega(G/K)$ of the form
$\mu(\theta)=\gamma(\theta)K$, as in Lemma \ref{lemmafirst}.
Since $a\in K$, we can write
$$(a.\mu)(\theta) = a\mu(\theta) = a\gamma(\theta)K
= a\gamma(\theta)a^{-1}K.$$
Then
\begin{align*}{}&\widetilde{(a\mu)}(\theta) = \sigma(a\gamma(\pi-\theta)a^{-1})K
=\sigma(a)\sigma(\gamma(\pi-\theta))\sigma(a^{-1})K
\\
{}& \ \ \ \ \ \ \ \ \ \ \ \ \ \ \ \  \ \ \ \ \ \  \ \  \  \ \  \ \ \  \ \ \ \ \ \ \ \ 
=a\sigma(\gamma(\pi-\theta))a^{-1}K=a\sigma(\gamma(\pi-\theta))K
=a\tilde{\mu}(\theta).\end{align*}
\end{proof}

We consider again the action of $A\times S^1$ on $\Omega(G)$ given by Equations
(\ref{conju}) and (\ref{ei}). 
We also recall (see Equation (\ref{taug})) that  $\tau$ is the
involutive automorphism  of $\Omega(G)$
 given by
 $$\tau(\gamma)(\theta) = \sigma(\gamma(-\theta)),$$
 $\theta\in S^1$ (see Equation (\ref{taug})). %L02/18
  The following proposition shows that the $A\times S^1$ action and
 the involution $\tau$ are compatible in the sense of Duistermaat \cite{Du}.
 \begin{propo}\label{propo}
We have
$$\tau((a,z).\gamma)=(a^{-1},z^{-1}).\tau(\gamma)$$
for any $\gamma\in \Omega(G)$ and any $(a,z)\in A\times S^1$.
\end{propo}

\begin{proof} 
We take  the $A$ and $S^1$ actions separately.
First, if $a\in A$ then we have $\sigma(a)=a^{-1}$, thus
$$\tau(a.\gamma)(\theta) =\sigma(a\gamma(-\theta)a^{-1})
= \sigma(a)\sigma(\gamma(-\theta))\sigma(a^{-1})
= a^{-1}\sigma(\gamma(-\theta))a
= (a^{-1}.\tau(\gamma))(\theta).$$
Second, if $z=e^{i\varphi}$, then
\begin{align*}{}&\tau(z.\gamma)(\theta) 
= \sigma(\gamma(-\theta+\varphi)\gamma(\varphi)^{-1})
= \sigma(\gamma(-\theta+\varphi))\sigma(\gamma(\varphi)^{-1})
\\{}& \ \ \ \ \ \ \ \ \ \  \  \ \  \ \ \  \ \ \ \  \ \ \ \  \ \  \  \  \ \  \  \ \ \ \  \ \ \  \ \ \ \  \ 
= \tau(\gamma)(\theta-\varphi)\tau(\gamma)(-\varphi)^{-1}
= (z^{-1}.\tau(\gamma))(\theta).\end{align*}

\end{proof}

We deduce immediately as follows.

\begin{coro}\label{coro}
The fixed point set
 $$\Omega(G)^\tau:=\{\gamma\in \Omega(G) \ : \ \sigma(\gamma(\theta))=\gamma(-\theta),
 \theta\in S^1\}$$
 is  invariant under the action of
 $A_2\times \bZ_2=\{(a,z)\in A\times S^1 \ : \ a^2=1, z=1 \ {\rm or} \ z=-1\}$.
 \end{coro}
 
 The following proposition makes the connection between the spaces
 $\Omega(G)^\tau$ and $\Omega(G/K)$. It is an equivariant version of a result whose origins go back to Bott and Samelson \cite{Bo-Sa} (see also \cite{Mi},
 \cite{Ko}). 
  
\begin{propo}\label{equiv}  There is a homotopy equivalence
between   $\Omega(G/K)$ and $\Omega(G)^\tau$ which is
equivariant with respect to the  $A_2\times \bZ_2$ actions defined
in Lemma \ref{commute} and Corollary \ref{coro}.
 \end{propo}
 
 \begin{proof} We use the idea of \cite[Proposition 3.1.3]{Ko} (see also
 \cite[Section 5]{Mi}).
 The homotopy equivalence is the map
 $F:\Omega(G)^\tau \to \Omega(G/K)$ given by
 $$F(\gamma):=\gamma|_{[0,\pi]}K.$$
 This map is well defined
 since if $\gamma$ is in $\Omega(G)^\tau$ then
 $\gamma(\pi)=\sigma(\gamma(\pi))$, thus
 $\gamma(\pi)\in K$ and consequently
 $\gamma(\pi)K=\gamma(0)K = eK.$
 To prove that $F$ is a homotopy equivalence, we note
 that we can identify $\Omega(G)^\tau$ with  the space of all paths 
 $\beta:[0,\pi]\to G$
  with $\beta(0)=e$ and $\beta(\pi)\in K$.
  The map $F$ is given by $\beta\mapsto \beta K$,
  for all paths $\beta$ as above.
  This is a principal bundle 
 whose fiber is the group  $\{\beta: [0,\pi]\to K\ : \ \beta(0)=e\}$. Since the latter space is contractible, 
  $F$ is a homotopy equivalence, as desired.
  
  It remains to show that $F$ is $A_2\times \bZ_2$ equivariant.
  Only the $\bZ_2$-equivariance is non-trivial.
Let us consider $\gamma\in \Omega(G)^\tau$ and verify that
$$F((-1).\gamma)=\widetilde{F(\gamma)}.$$
Here the loop $(-1).\gamma$ is given by
$$((-1).\gamma)(\theta) =\gamma(\theta+\pi)\gamma(\pi)^{-1},$$
for all $\theta\in S^1$, see Equation (\ref{ei}). Thus we have
$$F((-1).\gamma)(\theta) = \gamma(\theta+\pi)\gamma(\pi)^{-1}K
=\gamma(\theta+\pi)K,$$
since $\gamma(\pi)\in K$ (see above). On the other hand, for any $\theta\in S^1$ 
we have
$$\widetilde{F(\gamma)}(\theta) = \sigma(\gamma(\pi-\theta))K
=\gamma(\theta-\pi)K=\gamma(\theta+\pi)K.$$
  Here we have used that
  $\tau(\gamma)=\gamma$, which implies that
   $\sigma(\gamma(\pi-\theta))
=\gamma(\theta-\pi)$. 
   \end{proof}

Finally we can spell out the details of the proof of Corollary \ref{bottsam}:
 it follows from Theorem \ref{secmain} by using 
Proposition \ref{equiv} above (in the particular situation when $A=T$).

\section{Examples and counterexamples}\label{lasts}

\subsection{Examples}
The basic assumption of this paper is that the involutive automorphism
$\sigma$ of the
simply connected and compact Lie group $G$ 
satisfies $\sigma(t)=t^{-1}$ for all $t$ in a maximal torus $T\subset G$.
In other words, if $K$ denotes the fixed point set of $\sigma$, the
Riemannian symmetric pair $(G,K)$ is of maximal rank:
by this we mean that the rank of the symmetric space $G/K$ is equal to the rank of $G$ (not to be confused with the situation when the homogeneous space
$G/K$ has maximal rank, which means ${\rm rank} \ G={\rm rank } \ K$). Each Lie group $G$ as above has
essentially one such involution $\sigma$. 
In the following table we  describe $\sigma$  when $G$ is one of the classical simply connected compact Lie groups:
in each case it is sufficient to describe the automorphism $\theta:=d\sigma_e$ of $\g$.
We are using   \cite[Chapter X, Section 2, Subsection 3]{He}.

\begin{tabular}{|l|l|l|}
	\hline
$G$& $\g$ &  $\theta:=d\sigma_e$     \\
	\hline
$SU(n)$& $\s\u(n)$: $n\times n$ complex  skew-Hermitean\newline   & $\theta(X)=\overline{X}$   \\
 &  matrices $X$ & (complex conjugation)\\
	\hline
$Spin(2n)$& $\s\o(2n)$: $2n\times 2n$ real skew-symmetric &  $\theta(X)=I_{n,n}XI_{n,n}$,\\
 &  matrices $X$ & where $I_{n,n}:=
 \left(
\begin{array}{ccc}
-I_n & 0\\
0 & I_n
\end{array}
\right)$ \\
	\hline
$Spin(2n+1)$&  $\s\o(2n+1)$: $(2n+1)\times (2n+1)$ real  &  $\theta(X)=I_{n+1,n}X
I_{n+1,n}$,\\
 &skew-symmetric  matrices $X$ &where $I_{n+1,n}:=
 \left(
\begin{array}{ccc}
-I_{n+1} & 0\\
0 & I_{n}
\end{array}
\right)$  \\
	\hline
$Sp(n)$ &$\s\p (n)$: $X=	\left(
\begin{array}{ccc}
Z_{11} & Z_{12}\\
-\overline{Z}_{12} & Z_{22}
\end{array}
\right),$&  $\theta(X)=J_nXJ_n^{-1}$, \\
 & where $Z_{ij}$ are $n\times n$ complex matrices, &where 
 $J_{n}:=
 \left(
\begin{array}{cccc}
0 & I_n\\
-I_n & 0
\end{array}
\right)$  \\
 &  $Z_{11}$ and $Z_{22}$ skew-Hermitean,&  \\
 &  and $Z_{12}$ symmetric&  \\
\hline
\end{tabular}

\noindent The pairs $(G,\sigma)$ in the table above correspond to the symmetric spaces $G/K$ of
type  $AI$, $BDI$ (with $p=q$ or $p=q+1$), and
$CI$: for the meaning of these types, that is, for the classification of the irreducible Riemannian symmetric spaces, we refer the reader   
to \cite[Table V, p. 518]{He} or
\cite[Table 2, pp. 312-313]{Be}. For the exceptional Lie groups, one can also consult
the last two tables: the maximal rank types are $EI, EV, EVIII, FI,$ and $G$.

\subsection{Counterexamples}
In the remaining part of this section we will show that the hypothesis which says that
the pair 
$(G,K)$ is of maximal rank is essential for the two main results of the paper.  
The notations established in Remark 1 following Theorem \ref{firstmain} are in force here. 

Let us start with Theorem \ref{firstmain}. We  show that there exist 
 simply connected  compact Lie groups $G$ with an involution $\sigma$ and
a maximal torus $T\subset G$ such that $\Phi(\Omega(G)^\tau)$ is strictly
contained in $\Phi(\Omega(G))$. 
We first recall  that, in general,  the vertices of the  polyhedron
$\Phi(\Omega(G))$ in $\t \oplus \bR$  are
$\Phi(\gamma_\xi)=(\xi,\frac{1}{2}|\xi|^2)$, where $\xi $ is in the integral lattice 
$I$ of $T$
and $\gamma_\xi:S^1\to T$, $\gamma_\xi(\theta)=\exp(\theta\xi)$, for all 
$\theta\in S^1$, is the corresponding group homomorphism (see \cite[Section 1, Remark 2]{At-Pr}). 
Pick $\xi_0$ in the integer lattice $I$ such that $d\sigma_e(\xi_0)\neq -\xi_0$ (we will comment below on the existence of
such $\xi_0$). Let $\gamma_0:S^1\to T$, $\gamma_0(\theta)=\exp(\theta\xi_0)$
be the corresponding group homomorphism and
consider $\Phi(\gamma_0)=(\xi_0,\frac{1}{2}|\xi_0|^2)$.
Assume that there exists $\gamma \in \Omega(G)^\tau$ such that $\Phi(\gamma)
=\Phi(\gamma_0)$. 
Then $\Phi(\gamma)$ is on the paraboloid of equation
$E=\frac{1}{2}|p|^2$ in $\t \oplus \bR$, hence $\gamma$ must be a group homomorphism $S^1\to T$
(by \cite[Section 1, Remark 3]{At-Pr}). Thus $\gamma$ is of the form $\gamma(\theta)=\exp(\theta \xi)$, for all $\theta\in S^1$, where
$\xi\in I$. 
A simple calculation shows that the condition $\tau(\gamma)=\gamma$ implies $d\sigma_e(\xi)=-\xi$, thus $\xi\in \a$.
From $\Phi(\gamma)=\Phi(\gamma_0)$ we deduce 
$$(\xi,\frac{1}{2}|\xi|^2)=(\xi_0,\frac{1}{2}|\xi_0|^2)$$
thus $\xi=\xi_0$, which  contradicts $d\sigma_e(\xi_0)\neq -\xi_0$. 

One can easily find examples of symmetric spaces $G/K$ for which there exists
$\xi_0\in I$ with $d\sigma_e(\xi_0)\neq -\xi_0$. For example, one can take  
$$\bC P^{n-1}=SU(n)/S(U(1)\times U(n-1)).$$ This is a rank 1 symmetric space
(cf. e.g. \cite[Chapter X, Section 6, Table V]{He} or \cite[Example 6.6]{Mi}).
 Recall that the rank of a general symmetric space $G/K$ is equal to the 
 dimension of $\a$
(cf. e.g. \cite[Chapter V, Section 6]{He}, see also Remark 1 following Theorem \ref{firstmain}).
Thus, in the case at hand we have $\dim \a=1$.
We can  extend $\a$ to a maximal abelian subspace of ${\rm Lie}(SU(n))$,
call it $\t$, which is $d\sigma_e$ invariant 
and such that  
$$\a=\{x\in \t \ : \ d\sigma_e(x)=-x\}.$$
Put $T=\exp(\t)$, which is a maximal torus in $SU(n)$.
 It is clear that if $n\ge 3$, 
then $\dim \t=n-1$ is at least 2, and so 
not all integral elements of
$T$ are in $\a$. We note that the pair
$(SU(n), S(U(1)\times U(n-1))$ is far from being of maximal rank, as
${\rm rank } \ SU(n)=n-1$, whereas ${\rm rank} \  \bC P^{n-1}=1$.  An even more extreme example is given by the pair
$(Sp(n), U(n))$ (see \cite[Chapter X, Section 2, Subsection 3]{He}). This pair is 
of maximal rank: indeed, ${\rm rank} \ Sp(n)/U(n)= {\rm rank} \  Sp(n)=n$.
Hence there exists a torus $T\subset Sp(n)$ with
$\sigma(t)=t^{-1}$, for all $t\in T$: Theorem \ref{firstmain} applies in this situation.
 However,  we also have ${\rm rank} \ Sp(n)={\rm rank} \  U(n) =n$,
 thus there exists another maximal torus in $Sp(n)$, call it $T'$, such that
 $T'\subset U(n)$. This implies that $d\sigma_e(\xi)=\xi$, for all
 $\xi \in {\rm Lie}(T')$; thus $d\sigma_e(\xi)\neq -\xi$, unless $\xi=0$.

Let us now turn to Theorem \ref{secmain}.  This time we  show that there exist 
a simply connected compact Lie group $G$ with an involution $\sigma$ 
such that $\dim H^{2q}(\Omega(G);\bZ_2)\neq \dim H^q(\Omega(G)^\tau;\bZ_2)$,
for some $q\ge 0$. Indeed, let us consider again the  pair
$(SU(n), S(U(1)\times U(n-1)))$: the corresponding symmetric space is
$SU(n)/(S(U(1)\times U(n-1))=\bC P^{n-1}$ (see  above). The $\bZ_2$ Poincar\'e series of
$\Omega(SU(n))^\tau$ and $\Omega(\bC P^{n-1})$ are the same, being equal to
${(1+t)}{(1-t^{2n-2})^{-1}}$ (see  \cite[Section 6, Example 6.6]{Mi}).
The $\bZ_2$ Poincar\'e series of $\Omega(SU(n))$ is
$[(1-t^2)(1-t^4) \ldots (1-t^{2n-2})]^{-1}$ (cf. e.g. \cite[Equation (13.2)]{Bo-Sa}).
Thus if $n\ge 3$, then we have $\dim H^4(\Omega(SU(n));\bZ_2)=2$, whereas
$\dim H^2(\Omega(SU(n))^\tau;\bZ_2)=0$.

   \bibliographystyle{abbrv}

\end{document}